\documentstyle[11pt,newlfont,amscd,amssymb]{amsart}
\newcommand{\nc}{\newcommand}

\nc{\one}{\mbox{\bf 1}}
\nc{\invtensor}{\underset{\leftarrow}{\otimes}}
\nc{\rlarrows}{\begin{picture}(1,0.4)
                \put(0,-0.1){\makebox(1,0.2){$\leftarrow$}}
                \put(0,0.1){\makebox(1,0.2){$\ra$}}
              \end{picture}}
\nc{\rra}{\begin{picture}(1,0.4)
                \put(0,-0.1){\makebox(1,0.2){$\lra$}}
                \put(0,0.1){\makebox(1,0.2){$\lra$}}
              \end{picture}}
\nc{\Left}{\mathbf L}  % for derived
\nc{\Right}{\mathbf R} % functors
\nc{\gr}{\operatorname{gr}}
\nc{\Ho}{\operatorname{Ho}}
\nc{\alt}{\operatorname{alt}}
\nc{\Sym}{\operatorname{Sym}}
\nc{\sym}{\operatorname{sym}}
\nc{\id}{\operatorname{id}}
\nc{\Def}{\operatorname{Def}}
\nc{\Del}{\operatorname{Del}}
\nc{\Der}{\operatorname{Der}}
\nc{\im}{\operatorname{Im}}
\nc{\Ker}{\operatorname{Ker}}
\nc{\hfib}{\operatorname{h-fib}}

\nc{\coker}{\operatorname{Coker}}
\nc{\Col}{\operatorname{Col}}
\nc{\ter}{\operatorname{ter}}
\nc{\intl}{\operatorname{int}}
\nc{\out}{\operatorname{out}}
\nc{\val}{\operatorname{val}}

\nc{\Norm}{\operatorname{N}}
\nc{\Nor}{\operatorname{N}}
\nc{\Tor}{\operatorname{Tor}}
\nc{\res}{\operatorname{res}}
\nc{\Stab}{\operatorname{Stab}}
\nc{\Hom}{\operatorname{Hom}}
\nc{\uhom}{\CH\!o\!m}
\nc{\End}{\operatorname{End}}
\nc{\holim}{\operatorname{holim}}
\nc{\dirlim}{\underset{\rightarrow}{\lim}\,}
\nc{\invlim}{\underset{\leftarrow}{\lim}\,}
\nc{\CB}{\operatorname{\bf CB}}
\nc{\com}{\operatorname{co}}
\nc{\Tot}{\operatorname{Tot}}
\nc{\Th}{\operatorname{Th}}
\nc{\Cech}{\check{C}}
\nc{\Spec}{\operatorname{Spec}}
\nc{\Spf}{\operatorname{Spf}}
\nc{\MC}{\operatorname{MC}}
\nc{\U}{\operatorname{U}}
\nc{\Diff}{{\cal D}\mbox{\em iff}}
\nc{\Mor} {{\cal M}or}
\nc{\Ob}{\operatorname{Ob}}
\nc{\cone}{\widehat}
\nc{\Coder}{\operatorname{Coder}}
\nc{\pr}{\operatorname{pr}}
\nc{\diag}{\operatorname{diag}}
\nc{\tar}{\operatorname{tar}}

\nc{\SCAT}{\mathtt{SCAT}}
\nc{\Mod}{{\mathtt{mod}}}       
\nc{\Modf}{{\mathtt{modf}}}       
\nc{\Modg}{{\mathtt{modg}}}       
\nc{\Ab}{{\mathtt {Ab}}}          
\nc{\Alg}{{\mathtt {Alg}}} 
\nc{\Algf}{{\mathtt {Algf}}} 
\nc{\Algg}{{\mathtt {Algg}}} 
\nc{\Coalg}{{\mathtt {Coalg}}} 
         
\nc{\dgc}{{\mathtt{dgc}}}
\nc{\dgca}{{\mathtt{dgca}}}
\nc{\dgcu}{{\mathtt{dgcu}}}
\nc{\dgcuf}{{\mathtt{dgcuf}}}
\nc{\dgcf}{{\mathtt{dgcf}}}
\nc{\dgcg}{{\mathtt{dgcg}}}
\nc{\dgcc}{{\mathtt{dgccc}}}

\nc{\dgl}{{\mathtt{dglie}}}
\nc{\dgla}{{\mathtt{dgla}}}
\nc{\dglf}{{\mathtt{dglf}}}
\nc{\dglg}{{\mathtt{dglg}}}
\nc{\dga}{{\mathtt{dga}}}
\nc{\art}{{\mathtt {art}}}
\nc{\dgar}{{\mathtt {dgart}^{\leq 0}}}
\nc{\simpl}{\Delta^0\Ens}
\nc{\Coll}{{\mathtt{Coll}}}
\nc{\Kan}{{\mathtt {Kan}}}

\nc{\Grp}{{\mathtt {Grp}}}
\nc{\Cat}{{\mathtt {Cat}}}
\nc{\sGrp}{{\mathtt {sGrp}}}
\nc{\sCat}{{\mathtt {sCat}}}
\nc{\sFun}{{\mathtt {sFun}}}
\nc{\Ens}{{\mathtt {Ens}}}
\nc{\op}{{\operatorname{op}}}
\nc{\Op}{{\mathtt{Op}}}

\nc{\Lie}{{\mathtt{LIE}}}
\nc{\Com}{{\mathtt{COM}}}
\nc{\pa}{\partial}

\nc{\CA}{\cal A}
\nc{\CC}{\cal C}
\nc{\CDD}{\cal D}
\nc{\CF}{\cal F}
\nc{\CG}{\cal G}
\nc{\CH}{\cal H}
\nc{\CI}{\cal I}
\nc{\CJ}{\cal J}
\nc{\CL}{\cal L}
\nc{\CM}{\cal M}
\nc{\CN}{\cal N}
\nc{\CO}{\cal O}
\nc{\CP}{\cal P}
\nc{\CS}{\cal S}
\nc{\CT}{\cal T}
\nc{\CU}{\cal U}
\nc{\CW}{\cal W}
\nc{\fa}{\frak a}
\nc{\fg}{\frak g}
\nc{\fk}{\frak k}
\nc{\fh}{\frak h}
\nc{\fm}{\frak m}
\nc{\fn}{\frak n}
\nc{\fS}{\frak S}
\nc{\fI}{\frak I}
\nc{\fA}{\frak A}

\nc{\nen}{\newenvironment}
\nc{\ol}{\overline}
\nc{\ul}{\underline}
\nc{\lra}{\longrightarrow}
\nc{\lla}{\longleftarrow}
\nc{\Lra}{\Longrightarrow}
\nc{\Lla}{\Longleftarrow}
\nc{\Llra}{\Longleftrightarrow}
\nc{\hra}{\hookrightarrow}
\nc{\iso}{\overset{\sim}{\lra}}

\nc{\Thm}[1]{Theorem~\ref{#1}}
\nc{\Prop}[1]{Proposition~\ref{#1}}
\nc{\Lem}[1]{Lemma~\ref{#1}}
\nc{\Cor}[1]{Corollary~\ref{#1}}
\nc{\Conj}[1]{Conjecture~\ref{#1}}
\nc{\Claim}[1]{Claim~\ref{#1}}
\nc{\Defn}[1]{Definition~\ref{#1}}
\nc{\Exa}[1]{Example~\ref{#1}}
\nc{\Rem}[1]{Remark~\ref{#1}}
\nc{\Note}[1]{Note~\ref{#1}}

\nen{thm}[1]{\label{#1}{\bf Theorem.\ } \em}{}
\nen{prop}[1]{\label{#1}{\bf Proposition.\ } \em}{}
\nen{lem}[1]{\label{#1}{\bf Lemma.\ } \em}{}
\nen{klem}[1]{\label{#1}{\bf Key Lemma.\ } \em}{}
\nen{cor}[1]{\label{#1}{\bf Corollary.\ } \em}{}
\nen{conj}[1]{\label{#1}{\bf Conjecture.\ } \em}{}
\nen{claim}[1]{\label{#1}{\bf Claim.\ } \em}{}

\nen{defn}[1]{\label{#1}{\bf Definition.\ } }{}
\nen{exa}[1]{\label{#1}{\bf Example.\ } }{}

\nen{rem}[1]{\label{#1}{\em Remark.\ } }{}
\nen{note}[1]{\label{#1}{\em Note.\ } }{}
\nen{exer}[1]{\label{#1}{\em Exercise.\ } }{}

\begin{document}

%  top matter
\title[]{Deformations of homotopy algebras}
\author{Vladimir Hinich}
\address{Dept. of Mathematics, University of Haifa,
Mount Carmel, Haifa 31905 Israel}

%\thanks{}
\maketitle

%\tableofcontents
\section{Introduction}

\subsection{}
Let $k$ be a field of characteristic zero, $\CO$ be a dg operad over $k$
and let $A$ be an $\CO$-algebra. In this note we define formal deformations
of $A$, construct the deformation functor
$$\Def_A:\dgar(k)\to\simpl$$
from the category of  artinian local dg algebras (see~\ref{notation} for
the precise definition) to the category of simplicial sets. In the case
$\CO$ and $A$ are non-positively graded, we prove that $\Def_A$
is governed by the tangent Lie algebra $T_A$ defined in~\cite{haha}.
A very easy example~\ref{void} shows that the result does not hold
without this condition.

\subsection{}
``Classical'' formal deformation theory over a field of characteristic zero 
can be described as follows.

Let $k$ be a field of characteristic zero, $\art(k)$ be the category of
artinian local $k$-algebras with residue field $k$. Let $\CC$ be a category
cofibred over $\art(k)$. Equivalently, this means that a 2-functor
$$ \CC:\art(k)\to \Cat$$
is given, that is a collection of categories $\CC(R),\ R\in\art(k)$,
of functors $f^*:\CC(R)\to\CC(S)$ for each morphism $f:R\to S$ in $\art(k)$
and of isomorphisms $f^*g^*\iso (fg)^*$ satisfying the cocycle condition.

Finally, let an object $A\in\CC(k)$ be given. Then the deformation functor
$$\Def_A:\art(k)\to\Grp$$
assigns to each $R\in\art(k)$ the groupoid whose objects are isomorphisms
$\alpha:\pi^*(B)\to A$ where $\pi:R\to k$ is the natural map, and morphisms are
isomorphisms $B\to B'$ compatible with $\alpha$ and $\alpha'$.

If we are lucky enough (and this is usually the case), one can naturally
assign to $A$ a dg Lie algebra $\fg$ which  governs deformation of $A$
in the following sense: there is a natural equivalence of groupoids
$$ \Def_A\to\Del_{\fg}$$
as functors from $\art(k)$ to $\Grp$, where $\Del_{\fg}$ is the Deligne
groupoid of $\fg$ -- see~\cite{gm,ddg}.

As we argued in~\cite{uc}, the existence of dg Lie algebra $\fg$ 
is equivalent to the existence of formal dg moduli --- given by the
dg coalgebra $\CC(\fg)$. This also means that the functor $\Def_A$ can
be extended to a functor on the category $\dgar(k)$ of non-positively graded
dg commutative artinian algebras with residue field $k$ with values in
$\simpl$.

\subsection{}The category of algebras over a given dg operad $\CO$
is not just a category --- there exist weak equivalences, homotopies and
higher homotopies between the algebras. Therefore, the above
described approach can not produce a reasonable definition of formal
deformations of operad algebras.

Also, we want that deformations of quasi-isomorphic
algebras be equivalent, as well as deformations of algebras over
quasi-isomorphic operads, so we should use a sort of ``derived'' notion
of the deformation groupoid. 

\subsubsection{}
Morally, the picture should be the following.

For each $R\in\dgar(k)$ a $\infty$-category of $R\otimes\CO$-algebras
should be defined; denote it $\Alg^{\infty}(\CO,R)$. The collection of
$\Alg^{\infty}(\CO,R)$ should form an $\infty$-category cofibred 
over $\dgar(k)$.

Let now $A\in\Alg^{\infty}(\CO,k)$. Then the deformation functor
$$\Def_A:\dgar(k)\to\Grp^{\infty}$$
should be a ($\infty$-) functor to $\infty$-groupoids; its objects are
$\infty$-isomorphisms $\alpha:\pi_*(B)\to A$ and morphisms --- 
$\infty$-isomorphisms $B\to B'$ commuting with $\alpha$ and $\alpha'$.

\subsubsection{}
Since we do not know well what an $\infty$-category is and 
how to assign an $\infty$-category to the category of operad algebras,
we are looking for an appropriate substitute.

According to ~\cite{haha}, the category $\Alg(\CO,R)$ of 
$R\otimes\CO$-algebras admits a simplicial closed model category (SCMC) 
structure.

As a substitute to the $\infty$-category $\Alg^{\infty}(\CO,R)$, 
we suggest the simplicial category $\Alg^c_*(\CO,R)$ of cofibrant 
$R\otimes\CO$-algebras.

This allows us to define a deformation groupoid  as a functor
\begin{equation}
\Def_A:\dgar(k)\to\simpl
\label{def-fr-eq}
\end{equation}
to the simplicial sets.

Now, according to the general philosophy of deformation theory,  the
functor $\Def_A$ should be equivalent to the nerve of a certain
dg Lie algebra $\fg$ --- see~\cite{uc}, Sect.~8 ---  which is ``responsible'' 
for the deformations of $A$.  In~\cite{haha}, Sect.~7, we constructed a 
functorial tangent dg Lie algebra $T_A\in\dgl(k)$ as the Lie algebra 
$\Der(\widetilde{A},\widetilde{A})$ of derivations of a cofibrant
resolution $\widetilde{A}$ of $A$. The main result of this paper says that
$T_A$ governs the formal deformations of $A$ provided 
$\CO$ and $A$ belong to $C^{\leq 0}(k)$. Quite unexpectedly,
there is a very simple counter-example
showing that the last condition is necessary --- see Example~\ref{void}.

\subsection{Content of the Sections}
Throughout the paper we work a lot with simplicial categories and simplicial
groupoids. We collect in the Appendix the necessary information about
the subject. It is mostly well-known or easily imaginable; thus,
the closed model category structure on simplicial categories is a slight
generalization of a result of~\cite{dk}.

In Section~\ref{def-fr} we construct the deformation functor~ 
(\ref{def-fr-eq}).
To compare it with the nerve of the tangent dg Lie algebra, we provide
in Section~\ref{newnerve} a version of the nerve construction of
  ~\cite{uc}, Sect. 8, which assigns to a dg Lie algebra $\fg$ and to
a dg artinian algebra $R$ a simplicial groupoid.

Finally, in Section~\ref{final} we prove the main result. It follows easily 
from the model category structure of the category of simplicial categories.

\subsection{Notation}
\label{notation}

In what follows we use the following notations for different categories.

$\Ens,\Grp,\Cat$ are  the categories of sets, small groupoids and small 
categories respectively.

$\Delta$ is the category of ordered  sets $[n]=\{0,\ldots,n\}, \ n\geq 0$
and order-preserving maps. For a category $\CC$ we denote by $\Delta^0\CC$
the category of simplicial objects in $\CC$.

A simplicial category (and, in particular, a simplicial groupoid) will be
supposed to have a discrete set of objects, if it is not explicitly 
specified otherwise. The category of small simplicial categories is
defined $\sCat$ and that of simplicial groupoids $\sGrp$. 

For a fixed field $k$ of characteristic zero $\dgl(k)$ is the category 
of dg Lie algebras and $\dgar(k)$ is the category of non-positively graded
commutative artinian dg algebras with residue field $k$.

\subsection{Acknowledgement}This work was made during my stay at the 
Max-Planck Institut f\"ur Mathematik at Bonn. I express my gratitude
to the Institute for the hospitality.

\section{Homotopy algebras}
\label{def-fr}

In this Section we define  deformations of dg operad algebras.

Let $k$ be a fixed field of characteristic zero, $\CO\in\Op(C(k))$
be a dg operad over $k$, $A\in\Alg(\CO)$ be a $\CO$-algebra. 

Let $R$ be a commutative dg $k$-algebra. Componentwise tensoring by $R$  
defines a functor on $\Op(C(k))$ (its values can be equally considered
in $\Op(C(k))$ and in $\Op(\Mod(R))$). We define by $\Alg(\CO,R)$ the
category of $R\otimes\CO$-algebras. This category admits a CMC structure
--- see~\cite{haha} with quasi-isomorphisms as weak equivalences and
surjective maps as fibrations. We denote by $\Alg^c(\CO,R)$ the full
subcategory of cofibrant algebras and by $\CW^c(\CO,R)$ the category
of cofibrant $\CO$-$R$-algebras with weak equivalences as arrows.

The category $\Alg(\CO)$ admits also a simplicial structure so that
Quillen's axiom (SM7) is satisfied --- see~\cite{haha}, 4.8. The
simplicial structure is defined by the simplicial path functor
which assigns to an algebra $A\in\Alg(\CO)$ and to a finite simplicial set
$S\in\simpl$ the algebra $A^S=\Omega(S)\otimes A$ where $\Omega(S)$
denotes the dg commutative algebra of polynomial differential forms on $S$.

In the sequel we will add asterisk to denote that we consider the corresponding
simplicial category. Thus, for example,  $\CW^c_*(\CO,R)$ is the simplicial
category whose objects are cofibrant $R\otimes\CO$-algebras and whose
$n$-morphisms from $x$ to $y$ consist of quasi-isomorphisms
$x\to \Omega(\Delta^n)\otimes y$.

Recall~\cite{uc} that it is worthwhile to consider artinian local 
non-positively graded dg $k$-algebras $(R,\fm)\in\dgar(k)$ as bases
of formal deformations. 

\subsection{}
\begin{defn}{dha}Let $A\in\Alg(\CO)$. {\em Deformation functor}
 $$\Def_A:\dgar(k)\to\simpl$$
is defined by the formula
\begin{equation}
\Def_A(R)=\hfib_A\left(\CN(\CW^c_*(\CO,R))\to\CN(\CW^c_*(\CO,k))\right).
\label{def-def-functor}
\end{equation}
Here $\CN:\sCat\to\simpl$ is the simplicial nerve functor (see ~\ref{snerve})
and the homotopy fiber $\hfib$ being taken at a point 
$\widetilde{A}\in\CW^c_*(\CO,k)$
where $\widetilde{A}\to A$ is a cofibrant resolution of $A$.
\end{defn}

\subsubsection{}Recall (see~\cite{haha}) that for a $\CO$-algebra $A$ 
its tangent Lie algebra $T_A$ is defined as
$$ T_A=\Der(\widetilde{A},\widetilde{A}).$$

Now we are ready to formulate the main result of this note.
\subsubsection{}
\begin{thm}{main} Let $\CO$ be a dg operad over a field $k$ of characteristic
zero and let $A$ be an $\CO$-algebra. Suppose that both $\CO$
and $A$ are non-positively graded. Then the deformation functor
$\Def_A:\dgar(k)\to\simpl$ is equivalent to the nerve $\Sigma_{\fg}$
of the tangent dg Lie algebra $\fg:=T_A$. 
\end{thm}

\Thm{main} will be proven in Section~\ref{final}. 
Now we give an elementary example
which shows that the non-positivity condition is necessary.

\subsection{Example}
\label{void}
Let $\CO$ be the trivial operad $\CO(1)=k\cdot 1,\ \CO(i)=0$ for $i\not=0$.
$\CO$-algebras are just complexes and derivations are just all endomorphisms.

Let $A$ be the complex with zero differential with $A^i=k$ for all
$i\in\Bbb{Z}$. 

Then $T_A=\Hom(A,A)$ is a complex with zero differential; an element
$f\in(\fm\otimes T_A)^1$ satisfies the Maurer-Cartan equation iff $f^2=0$.
For instance, put $R=k[\epsilon]/\epsilon^2$. Then any element $f$ of degree
one is Maurer-Cartan. The corresponding to $f$ complex of $R$-modules
is $R\otimes A$ as a graded $R$-module, and has $f$ as the differential.
Suppose  that all components of $f$ are non-zero. Then $(R\otimes A,f)$ is contractible
and of course can not be thought of being a deformation of $A$. 
\section{Simplicial Deligne groupoid}
\label{newnerve}

\subsection{Definition}
Let $k$ be a field of characteristic zero and  $\fg\in\dgl(k)$ be a 
nilpotent dg Lie $k$-algebra. 
In this Section we construct a simplicial groupoid 
$\Gamma(\fg)=\{\Gamma_n(\fg)\}$ whose nerve (see~\ref{snerve}) is naturally 
homotopically equivalent to the nerve $\Sigma(\fg)$.

The construction is a generalization (and a simplification) of the one we used 
in~\cite{uc}, 9.7.6.

Recall (see~\cite{ddg},\cite{uc}, 8.1.1) that the nerve $\Sigma(\fg)$ 
of the nilpotent dg Lie algebra $\fg$ is defined as
\begin{equation}
\Sigma_n(\fg)=\MC(\Omega_n\otimes\fg),
\label{nerve}
\end{equation}
$\Omega_n$ being the algebra of polynomial differential forms on the standard
$n$-simplex.

Following~\cite{uc}, Sect. 8, define a simplicial group $G=G(\fg)$
by the formula
\begin{equation}
G_n=\exp(\Omega_n\otimes\fg)^0.
\label{G.}
\end{equation}
Here $\Omega_n\otimes\fg$ is a nilpotent dg Lie algebra, so its zero
component is an honest nilpotent Lie algebra, and therefore its exponent 
makes sense.

Define a simplicial groupoid $\Gamma:=\Gamma(\fg)$ (we will call it 
{\em simplicial Deligne groupoid} since its zero component is the 
conventional Deligne groupoid~\cite{gm}) as follows.

$\Ob\Gamma=\MC(\fg);$

$ \Hom_{\Gamma}(x,y)_n=\{g\in G_n|g(x)=y\}.$

It is useful to have in mind the following easy

\subsubsection{}
\begin{lem}{G-contractible}
The simplicial group $G(\fg)$ is always contractible. 
\end{lem}
\begin{pf}
As a simplicial set, $G$ is isomorphic to the simplicial vector space
$$ n\mapsto (\Omega_n\otimes\fg)^0.$$
The latter is a direct sum of simplicial vector spaces of form 
$\Omega_{\bullet}^p$ (each one $\dim\fg^{-p}$ times) which are all
contractible --- see~\cite{l}, p.~44.
\end{pf}

\subsection{Equivalence}

Recall that any simplicial category (and more generally, any 
$\CC\in\Delta^0\Cat$)
defines a bisimplicial set whose diagonal is called {\em the nerve}
of $\CC$, denoted by $\CN(\CC)$ --- see~\ref{snerve}.

\subsubsection{}
\begin{prop}{eq.defns}
The nerve $\Sigma(\fg)$ of a nilpotent dg Lie algebra is naturally
homotopically equivalent to $\CN(\Gamma(\fg))$.
\end{prop}
\begin{pf}Define $\Gamma'\in\Delta^0\Grp$ (a simplicial groupoid in the wide
sense) by the following formulas.

$ \Ob\Gamma'_n=\MC(\Omega_n\otimes\fg);$

$ \Hom_{\Gamma'}(x,y)_n=\{g\in G_n(\fg)|g(x)=y\}.$

One has a natural fully faithful embedding 
$\Gamma(\fg)\to\Gamma'$. According to~\cite{uc}, 
8.2.5, the map $\Gamma_n(\fg)\to\Gamma'_n$ is an equivalence of groupoids 
for each $n$. This implies that the induced map of the nerves 
$$\CN(\Gamma)\to\CN(\Gamma')$$
is a homotopy equivalence.

Now we shall compare the nerve $\CN(\Gamma')$ to $\Sigma(\fg)$. Look at
$\Gamma'$ as at a bisimplicial set. One has
$$ \Gamma'_{pq}=\Sigma_p(\fg)\times G_p(\fg)^q.$$

This means that the simplicial set $\Gamma'_{\bullet q}$ is equal to 
$\Sigma(\fg)\times G(\fg)^q$. 

The simplicial set $G(\fg)$ is contractible by~\Lem{G-contractible}. 
Therefore, $\Gamma'_{\bullet q}$ is canonically homotopy equivalent 
to $\Sigma(\fg)$. This implies that the nerve $\CN(\Gamma')$ is homotopy 
equivalent to $\Sigma(\fg)$.
\end{pf}

\subsubsection{}
\begin{rem}{gen}\Prop{eq.defns} generalizes the claim used 
in the proof of 9.7.6 of ~\cite{uc}. 
\end{rem}

\subsubsection{}
Let now $\fg\in\dgl(k)$. Following the well-known pattern, we define
the functor 
$$\Gamma_{\fg}:\dgar(k)\to\sGrp$$
by the formula 
$$\Gamma_{\fg}(R)=\Gamma(\fm\otimes\fg)$$
for $(R,\fm)\in\dgar(k)$.

The functor $\Gamma_{\fg}$ is also called the simplicial Deligne groupoid.

\subsection{Properties}

We wish to deduce now some properties of the simplicial Deligne groupoid 
functor which are similar to the properties of the nerve $\Sigma(\fg)$ ---
see~\cite{uc}, Sect.~8. 

In what follows we use the closed model category (CMC) structure
on the category $\sCat$ --- see~\ref{scat-cmc-ss}. 

\subsubsection{}
\begin{prop}{gamma-fib}
Let $f:\fg\to \fh$ be surjective (resp., a surjective quasi-isomorphism).
Then for each $(R,\fm)\in\dgar(k)$ the map
$$f:\Gamma_{\fg}(R)\to \Gamma_{\fh}(R)$$
is a fibration (resp., an acyclic fibration) in $\sCat$.
\end{prop}
\begin{pf}
Note first of all that the similar claim holds for the nerve functor: 
according to~\cite{uc}, Prop. 7.2.1, the map 
$f:\Sigma_{\fg}(R)\to\Sigma_{\fh}(R)$ is a fibration
(resp., acyclic fibration) provided $f$ is a surjection (resp., a surjective
quasi-isomorphism). This implies that the map
$$f:\Gamma_{\fg}(R)\to \Gamma_{\fh}(R)$$
satisfiest the property (1) of fibrations (resp., of acyclic fibrations) --- 
see~\ref{scat-fib},~\ref{scat-af}.

Let us check the property (2). It claims that for any 
$x,y\in\Ob\Gamma_{\fg}(R)$ the map of simplicial sets
$$ f:\uhom_{\fg}(x,y)\to\uhom_{\fh}(fx,fy)$$
is a Kan fibration (resp., acyclic Kan fibration) --- here
we write $\uhom_{\fg}(\_,\_)$ instead of $\uhom_{\Gamma_{\fg}(R)}(\_,\_)$.

Denote $G=G(\fg),\ H=G(\fh)$ the simplicial groups corresponding to 
$\fg$, $\fh$ as in the formula~\ref{G.}.

A map from a simplicial set $K$ to $\uhom_{\fg}(x,y)$ is given by an
element $g\in G(K)=\Hom(K,G)$ satisfying the condition $g(x)=y$.

Let a commutative diagram in $\simpl$
\begin{center}
$$\begin{CD}
K@>>> \uhom_{\fg}(x,y) \\
@V{\alpha}VV    @V{f}VV \\
L@>>> \uhom_{\fh}(fx,fy) \\
\end{CD}$$
\end{center}
be given with $\alpha:K\to L$ being a cofibration of finite simplicial sets.
We suppose also that either $\alpha$ or $f$ is a weak equivalence.
Our aim is to find a map $L\to \uhom_{\fg}(x,y)$ commuting with
the above diagram. Thus, we are given with a compatible pair of
elements $g\in G(K),\ h\in H(L)$ satisfying the groperty
$$ g(x)=y; \ h(fx)=fy.$$
Our aim is to lift this pair to an element $\widetilde{g}\in G(L)$
satisfying the property $\widetilde{g}(x)=y$.

We will do this in two steps. First of all, since $f$ is surjective, 
the induced map of simplicial groups $f: G\to H$ is surjective, and, 
therefore, fibrant.
Furthermore, since both $G$ and $H$ are contractible by~\Lem{G-contractible}, 
the map $f:G\to H$ is actually an acyclic fibration, and therefore  the pair 
of compatible elements $g\in G(K),\ h\in H(L)$ lifts to an element 
$g'\in G(L)$. We can not, unfortunately, be sure that $g'(x)=y$. This is why
we need the second step which will correct $g'$ to satisfy this property.

Suppose $g'(x)=y'\in\MC(\Omega(L)\otimes\fm\otimes\fg)$. The elements $y$
and $y'$ of $\MC(\Omega(L)\otimes\fm\otimes\fg)$ have the same images
in both $\MC(\Omega(K)\otimes\fm\otimes\fg)$ and 
$\MC(\Omega(L)\otimes\fm\otimes\fh)$. Now, the commutative diagram
\begin{center}
$$\begin{CD}
\Omega(L)\otimes\fg@>>> \Omega(K)\otimes\fg\\
@VVV    @VVV \\
\Omega(L)\otimes\fh@>>> \Omega(K)\otimes\fh\\
\end{CD}$$
\end{center}
induces an acyclic fibration 
$$ p:\fg_1:=\Omega(L)\otimes\fg\to\Omega(K)\otimes\fg\times
_{\Omega(K)\otimes\fh}\Omega(L)\otimes\fh=:\fg_2$$
of dg Lie algebras. Then
the map $\Sigma_p:\Sigma_{\fg_1}(R)\to \Sigma_{\fg_2}(R)$ is an acyclic
fibration.

 Now, we have two elements $y,y'\in\MC(\fm\otimes\fg_1)$ 
satisfying $p(y)=p(y')\in\MC(\fm\otimes\fg_2)$. Therefore, there
exists an element $z\in \Sigma_{\fg_1}(R)_1$ such that $d_0z=y,d_1z=y'$
and $p(z)=s_0(p(y))$. Using the explicit description of $\Sigma_{\fg}(R)_1$
in~\cite{uc},  8.2.3, one obtains and element $\gamma\in\exp(\fm\otimes\fg_1)$
satisfying $p(\gamma)=1\in\exp(\fm\otimes\fg_2);\ \gamma(y')=y$.

Then one immediately sees that the element $\widetilde{g}=\gamma g'$ is the one
we need.
\end{pf}

\subsubsection{}
\begin{cor}{}
1. For any $\fg\in\dgl(k), R\in\dgar(k), x,y\in\Ob\Gamma_{\fg}(R)$ the 
simplicial set $\uhom(x,y)$ is fibrant.

2. Any quasi-isomorphism $f:\fg\to\fh$ induces a weak equivalence
$$ f:\Gamma_{\fg}(R)\to\Gamma_{\fh}(R)$$
for each $R\in\dgar(k)$.
\end{cor}
\begin{pf}
1. Take $\fh=0$ in~\Prop{gamma-fib}.

2. The category $\dgl(k)$ admits a CMC structure with surjections as 
fibrations and quasi-isomorphisms as weak equivalences --- see~\cite{haha},
Sect.~4. Using this, present $f=p\circ i$ as a composition of an acyclic
fibration $p$ and an acyclic cofibration $i$. Any acyclic cofibration in 
$\dgl(k)$ is left invertible: $q\circ i=\id$. The map $q$ is obviously an
acyclic fibration. Then by~\Prop{gamma-fib} the map 
$f:\Gamma_{\fg}(R)\to\Gamma_{\fh}(R)$ is a weak equivalence.
\end{pf}

\section{Final}
\label{final}

\subsection{}
We start with an observation explaining the connection between $T_A$ and
the formal deformations of $A$. Let $B$ be a cofibrant $R\otimes\CO$-algebra
with $(R,\fm)\in\dgar(k)$.
Denote $A=k\otimes_RB$. The algebra $B$ is isomorphic, as a graded
$\CO$-algebra, to $R\otimes A$. Choose a graded isomorphism
$$ \theta: B\to R\otimes A$$
and put
$$ z=\theta\circ d_B\circ\theta^{-1}-1\otimes d_A$$
where $d_B$ (resp., $d_A$) is the differential in $B$ (resp., in $A$).

Then $z$ is a degree one derivation in $\fm\otimes T_A$ satisfying the
Maurer-Cartan equation. A different choice of isomorphism $\theta$
gives rise to a Maurer-Cartan element $z'\in\fm\otimes T_A$ equivalent to
$z$: there exists $g\in\exp(\fm\otimes T_A)^0$ such that $z'=g(z)$.

In what follows we will use a (non-unique) presentation of a  
$R\otimes\CO$-algebra $B$ by an element $z\in\MC(\fm\otimes T_{k\otimes_RB})$.

\subsection{Proof of the Theorem}

To simplify the notation, denote $\CW=\CW^c_*(\CO,R),\ol{\CW}=\CW^c_*(\CO,k)$.

\subsubsection{}
\begin{lem}{isafibration}
The natural map $\pi:W\to\ol{W}$ is a fibration in $\sCat$.
\end{lem}
\begin{pf}
1. Let us prove the condition (1) of Definition~\ref{scat-fib}. 

It means the following.
Let $f:A\to B$ be a quasi-isomorphism of cofibrant $\CO$-algebras over $k$.
Let one of two elements $a\in\MC(\fm\otimes T_A)$ or $b\in\MC(\fm\otimes T_B)$
be given. We have to check that there exists a choice of the second element
and a map 
$$g: (R\otimes A,d+a)\to(R\otimes B,d+b)$$
of $R\otimes\CO$-algebras which lifts $f:A\to B$.

We can consider separately the cases when $f$ is an acyclic fibration or
an acyclic cofibration.

In both cases we will be looking for the map $g$ in the form 
$$g=\gamma^{-1}_B\circ (\id_R\otimes f)\circ\gamma_A$$
where $\gamma_A\in\exp(\fm\otimes T_A)^0$ and similarly for $\gamma_B$.
A map $g$ as above should commute with the differentials $d+a$ and $d+b$.
This amounts to the condition
$$ f_*(\gamma_A(a))=f^*(\gamma_B(b)),$$
where the natural maps 
$$ T_A\overset{f_*}{\lra}\Der_f(A,B)\overset{f^*}{\lla}T_B$$
are defined as in~\cite{haha}, 8.1.

Recall that we are assuming that $f$ is either acyclic cofibration or an 
acyclic fibration.

In  both cases there exists a commutative square

\begin{center}
$$\begin{CD}
T_f@>{\alpha}>> T_A \\
@V{\beta}VV    @V{f_*}VV \\
T_B@>{f^*}>> \Der_f(A,B)\\
\end{CD}$$
\end{center}
where $T_f$ is a dg Lie algebra and $\alpha,\beta$ are Lie algebra 
quasi-isomorphisms --- see~\cite{haha}, 8.2, 8.3. 
The maps $\alpha,\beta$ induce bijections
$$ \pi_0(\Sigma_{T_A}(R))\lla\pi_0(\Sigma_{T_f}(R))
\lra\pi_0(\Sigma_{T_A}(R))$$
which prove the assertion.

2. Let us check the condition (2) of~\ref{scat-fib}.
Let $\widetilde{A},\widetilde{B}\in\CW$ and let $A=k\otimes_R\widetilde{A},
B=k\otimes_R\widetilde{B}$. We have to check that the map
\begin{equation}
\uhom(\widetilde{A},\widetilde{B})\to\uhom(A,B)
\label{red}
\end{equation}
is a Kan fibration.
But this results from the SCMC structure on $\Alg(R\otimes\CO)$. In fact,
$\widetilde{A}$ is cofibrant, and the reduction map $\widetilde{B}\to B$
can be considered as a fibration in $\Alg(R\otimes\CO)$. Therefore,
the map
\begin{equation}
\uhom(\widetilde{A},\widetilde{B})\to\uhom(\widetilde{A},B)
\label{red2}
\end{equation}
is a Kan fibration. But the maps~(\ref{red2}) and (\ref{red}) coincide,
so the condition (2)  of~\ref{scat-fib} is proven.
\end{pf}

\subsubsection{}
Fix now a cofibrant $\CO$-algebra $A$ and denote $\fg=T_A$.
Fix $(R,\fm)\in\dgar(k)$.

Define a map of simplicial categories
$$ \alpha:\Gamma_{\fg}(R)\to W$$
as follows. Let $z\in\MC(\fm\otimes\fg)=\Ob\Gamma_{\fg}(R)$. Put
$$\alpha(z)=(R\otimes A, 1\otimes d+z).$$
Now, any element $g\in G_n=\exp(\Omega_n\otimes\fm\otimes\fg)^0$
defines a graded automorphism of $\Omega_n\otimes R\otimes A$. This
obviously defines an isomorphism of $R\otimes\CO$-algebras
$$ (R\otimes A, 1\otimes d+z)\lra (R\otimes A, 1\otimes d+g(z)).$$

The map $\alpha:\Gamma_{\fg}(R)\to W$ identifies $\Gamma_{\fg}(R)$
with the fibre of $\pi:W\to\ol{W}$ at $A$. Since $\pi$ is
a fibration, this is weakly equivalent to the homotopy fibre of $\pi$.

Now Theorem follows since the nerve functor $\CN$ preserves 
fibrations and weak equivalences by~\ref{nerve-exactness}.

\subsection{Concluding remarks}
\subsubsection{}
The formula~(\ref{def-def-functor}) defines a deformation functor for
operad algebras not necessarily concentrated in non-positive degrees.

The definition seems to be correct also in this case, in spite of
the fact that $T_A$ does not govern deformations of $A$ in this case.
It seems that one should find in this case a more clever way to define
the tangent Lie algebra. It might consist of ``tame''
derivations of a cofibrant resolution $\widetilde{A}$ of $A$, the ones 
which behave well with respect to a filtration on $\widetilde{A}$.

\subsubsection{}
One can define ``hard'' formal deformations of a $\CO$-algebra $A$ which deform
not only the algebra $A$ itself but also the base operad $\CO$. 
Then the universal deformation of $(\CO,A)$ would provide the
tangent Lie algebra $T_A$ with a canonical extra structure.

\section{Appendix: simplicial categories}

In this Section we present a more or less standard information about
simplicial categories. It includes the description~\ref{scat-cmc} of a CMC
structure on the category $\sCat$ of small simplicial categories.
This structure is a slight generalization of the one described in~\cite{dk}.

\subsection{Weak equivalences and fibrations in $\sCat$ }
\label{scat-cmc-ss}

Here we define a closed model category structure on the category $\sCat$
of simplicial categories.

\subsubsection{(Co)limits}
\label{lim}
The category $\sCat$ admits arbitrary limits and colimits.  
Inverse limits in $\sCat$ are induced by inverse limits in $\Ens$ in the
obvious sense. 

The existence of inductive limits in $\sCat$ follows by a general abstract 
nonsense from the existence of inductive limits in $\Ens$.
Note that the functor $\sCat\to\Ens$ assigning to each simplicial category
the set of its objects, commutes with inductive limits. The set of morphisms
of an inductive limit is freely generated by the morphisms of all categories
involved, modulo an obvious equivalence relation.

Note that the existence of direct limits in $\sCat$ allows one to mimic the
procedure of ``adding variables''. We will single out the following cases.

{\em Adding an object.} Given $\CC\in\sCat$, denote $\CC\langle *\rangle$
the coproduct of $\CC$ with the trivial one-object category $*$.

{\em Adding an ingoing arrow.} Given $\CC\in\sCat,\ x\in\Ob\CC$, one defines
$\CC\langle *\to x\rangle$ with the set of objects $\Ob\CC\coprod\{*\}$
and the set of morphisms freely generated by $\Mor\CC$ and by the map 
$*\to x$.

{\em Adding an outgoing arrow.} The category $\CC\langle x\to *\rangle$ is
defined similarly to the above.

{\em Adding maps between objects} Given $\CC\in\sCat$, $x,y\in\Ob\CC$
and a map $\alpha:\uhom_{\CC}(x,y)\to H$ of simplicial sets, the simplicial
category $\CC\langle x,y;\alpha\rangle$ has the same objects as $\CC$. Its
set of morphisms is freely generated by $\Mor\CC$ and by $H$.

\subsubsection{}
Define the functor 
$$\pi_0:\sCat\to\Cat$$
as follows. For $\CC\in\sCat$ the category $\pi_0(\CC)$ has the same objects
as $\CC$. For $x,y\in\Ob\pi_0(\CC)$

$$\Hom_{\pi_0(\CC)}(x,y)=\pi_0(\Hom_{\CC}(x,y)).$$

\subsubsection{}
\begin{defn}{scat-we}
A map $f:\CC\to\CDD$ in $\sCat$ is called a weak 
equivalence if the following properties are satisfied.

(1) The map $\CN(\pi_0(f))$ is a weak equivalence of simplicial sets.

(2) For all $x,x'\in\Ob\CC$ the map $f:\uhom(x,x')\to\uhom(fx,fx')$
is a weak equivalence.
\end{defn}

\subsubsection{}
\begin{defn}{scat-fib}A map $f:\CC\to\CDD$ in $\sCat$ is called a fibration
if it satisfies the following properties

(1) the right lifting property (RLP) with respect to ``adding an ingoing
or an outgoing arrow''     
$$\CA\to\CA\langle *\to x\rangle,\ \CA\to\CA\langle x\to *\rangle$$
 (see~\ref{lim}).

(2) For all $x,x'\in\Ob\CC$ the map $f:\uhom(x,x')\to\uhom(fx,fx')$
is a Kan fibration. This is equivalent the the RLP with respect to
all maps $\CA\to\CA\langle x,y;\alpha\rangle$ where $\alpha$ is an
acyclic fibration (see~\ref{lim}).
\end{defn}

\subsubsection{}
\begin{thm}{scat-cmc}The category $\sCat$ admits a CMC structure with
weak equivalences described in~\ref{scat-we} and fibrations as 
in~\ref{scat-fib}.
\end{thm}

\subsubsection{}
An explicit description of different classes of morphisms in $\sCat$ is given
 below. The proof of the Theorem is standard. It is given
in~\ref{pf-scat-cmc}. 

\subsubsection{}
A map $f:\CC\to\CDD$ in $\sCat$ is called an 
acyclic fibration if it is simultaneously a weak equivalence and a fibration.

\begin{lem}{scat-af}
A map $f:\CC\to\CDD$ is an acyclic fibration iff the following conditions
are satisfied.

(1) the map $\Ob f:\Ob\CC\to\Ob\CDD$ is surjective. In other words, 
$f$ satisfies the RLP with respect to ``adding an object map'' 
$\CA\to\CA\langle *\rangle$.

(2) For all $x,x'\in\Ob\CC$ the map $f:\uhom(x,x')\to\uhom(fx,fx')$
is an acyclic Kan fibration.
\end{lem}
\begin{pf}If $f$ satisfies (1), (2), it is clearly an acyclic fibration.
In the other direction, suppose $f$ is an acyclic fibration. Then the 
property (2) is clear. We have only to check that $\Ob f$ is surjective.
Since $f$ satisfies the RLP with respect to ingoing and outgoing arrows,
$\CDD$ is a disjoint union of the full subcategories, defined by
the image of $\Ob f$ and by its complement. Since $\pi_0(f)$ is a weak 
equivalence, it induces a bijection of the connected components of $\CC$ and 
$\CDD$ and this proves the claim.  
\end{pf}

\subsubsection{}
\label{sc}
A map $f:\CC\to\CDD$ will be called {\em a standard cofibration} if there is 
a collection of maps $f_i:\CC_i\to\CC_{i+1},\ i\in\Bbb{N}$ such that
$\CC=\CC_0,\ \CDD=\dirlim\CC_i$, and each $f_i$ is a coproduct of maps
of one of the following two types:

(1) Adding an object $\CC_i\to\CC_i\langle *\rangle$;

(2) Adding maps between objects $\CC_i\to\CC_i\langle x,y;\alpha\rangle$
 with $\alpha$ injective.

By~\ref{scat-af}, standard cofibrations satisfy the LLP with respect to all
acyclic fibrations.

\subsubsection{}
\label{sac}
A map $f:\CC\to\CDD$ will be called {\em a standard acyclic cofibration} 
if there is a collection of maps $f_i:\CC_i\to\CC_{i+1},\ i\in\Bbb{N}$ 
such that $\CC=\CC_0,\ \CDD=\dirlim\CC_i$, and each $f_i$ is a coproduct of 
maps of one of the following three types:

(1+) Adding an ingoing arrow $\CC_i\to\CC_i\langle *\to x\rangle$;

(1--) Adding an outgoing arrow $\CC_i\to\CC_i\langle x\to *\rangle$;

(2) Adding maps between objects $\CC_i\to\CC_i\langle x,y;\alpha\rangle$
 with $\alpha$ acyclic cofibration.

By~\ref{scat-fib}, standard acyclic cofibrations satisfy the LLP with respect
to all fibrations.

\subsubsection{}
The following description of cofibrations and of acyclic cofibrations
results from the proof of~\Thm{scat-cmc}.

\begin{cor}{c-and-ac}1. Any cofibration in $\sCat$ is a retract of a standard
cofibration.

2. Any acyclic cofibration in $\sCat$ is a retract of a standard acyclic
cofibration.
\end{cor}

\subsection{Proof of~\Thm{scat-cmc}}
\label{pf-scat-cmc}

The axioms (CM 1), (CM 2), (CM 3), (CM 4)(ii) are immediately verified.

(CM 5)(ii) Let $f:X\to Y$ be a map in $\sCat$.
Adding objects to $X$, we can ensure that 
the map $f:\Ob(X)\to\Ob(Y)$ is surjective.
Then, adding maps between objects, we can decompose $f$ into a
standard cofibration followed by an acyclic fibration.

This implies, in particular, that any cofibration is a retract of a standard
cofibration.

To check the axiom (CM 5)(i) we need the following
\subsubsection{}
\begin{lem}{sac-is-ac}
Standard acyclic cofibrations are acyclic cofibrations.
\end{lem}
\begin{pf}It is enough to prove that a map $\CC\to\CDD$ is a weak equivalence
when $\CDD$ is obtained from $\CC$ by one of the following ways.

(1) adding a number of ingoing arrows;

(2) adding a number of outgoing arrows;

(3) adding (simultaneously) maps between objects $x_i$ and $y_i$ along
acyclic cofibrations $\alpha_i:\uhom(x_i,y_i)\to H_i$.

In the first two cases the map $\CC\to \CDD$ is easily split by an acyclic
fibration.

The shortest way to get the result in the case (3) is to use Proposition 7.2
of~\cite{dk} which claims the existence of CMC structure on the category
of simplicial categories having a fixed set of objects. 

\end{pf}

(CM 5)(i)  Let $f:X\to Y$ be a map in $\sCat$. Adding ingoing and outgoing
arrows to $X$, we can ensure that the image of $\Ob(X)$ under $f$ consists
of a number of connected components of $\Ob(Y)$. From now on we can
suppose, without loss of generality, that $f$ is surjective on objects.
Then for a decomposition $f=p\circ i$ it is enough to check that $p$
satisfies condition (2) of ~\ref{scat-fib} to ensure $p$ is a fibration.

Now applying step by step the procedure of adding maps between objects 
$\CC\to\CC\langle x,y;\alpha\rangle$ along acyclic cofibrations $\alpha$,
we can construct a decomposition $f=p\circ i$ with $p$ fibration and $i$
a standard acyclic cofibration. According to~\Lem{sac-is-ac}, $i$ is an acyclic
cofibration.

Now, applying the proof of (CM 5)(i) to any acyclic cofibration $f$, we deduce
that $f$ is a retract of a standard acyclic cofibration.

(CM 4)(i) By definition, any standard acyclic cofibration satisfies LLP 
with respect to all fibrations. Any acyclic fibration is a retract
of a standard acyclic fibration, and therefore satisfies as well LLP
with respect to all fibrations.  

Theorem is proven.

\subsection{Simplicial nerve}
\label{snerve}
\subsubsection{}
\label{snerve}
In what follows we identify $\Cat$ with the full subcategory of
$\simpl$. Then every simplicial category (and even every $\CC\in\Delta^0\Cat$)
can be seen as a bisimplicial set; its diagonal will be called {\em
the nerve} of $\CC$ and will be denoted $\CN(\CC)$. If $\CC$ is a ``usual''
category, $\CN(\CC)$ is its ``usual'' nerve. 

The functor $\CN:\sCat\to\simpl$ admits a left adjoint functor
$$\SCAT:\simpl\to\sCat$$
defined by the properties

\begin{itemize}
\item{$\Ob\SCAT(\Delta^n)=[n]=\{0,\ldots,n\};$}
\item{$\Mor\SCAT(\Delta^n)\text{ is freely generated by }a_i\in\uhom_n(i-1,i),
\ i=1,\ldots,n;$}
\item{ $\SCAT$ commutes with arbitrary colimits.}
\end{itemize}

\subsubsection{}
\begin{prop}{nerve-exactness}
The nerve functor $\CN:\sCat\to\simpl$ preserves weak equivalences, fibrations
and cofibrations.
\end{prop}
\begin{pf}
1. To check that $\CN$ preserves the fibrations, it is enough to prove
that the adjoint functor $\SCAT$ preserves acyclic cofibrations. For this
we have to check that $\SCAT$ transforms any map $\Lambda^n_i\to\Delta^n$
to an acyclic fibration. This is an easy exercise (one should consider the 
cases $n=1$ and $n>1$ separately). Note that the same reasoning (even
easier!) proves that $\CN$ preserves acyclic fibrations --- this is because
$\SCAT$ preserves cofibrations.

2. It is clear that $\CN(f)$ is a weak equivalence provided $f$ is a weak
equivalence {\em bijective on objects}. To prove the general claim,
we present $f$ as a composition of an acyclic fibration with an acyclic 
cofibration and therefore reduce the problem to the case $f$ is an 
acyclic cofibration. Using~\ref{c-and-ac}, we can suppose that $f$ is
of one of the types (1+), (1--), (2) of~\ref{sac}. The type (2) does not
change the set of objects, so we have nothing to prove. The maps of types
(1+), (1--) split, and the splitting map is an acyclic fibration. This
proves the claim.

3. The claim about cofibrations is obvious.
\end{pf}

\end{document}